\documentclass[12pt,a4paper]{article}
\usepackage{amsmath, amssymb, theorem, latexsym, epsfig}

\usepackage{graphics}

\usepackage{epstopdf}

\newtheorem{theorem}{Theorem}[section]
\newtheorem{lemma}[theorem]{Lemma}
\newtheorem{proposition}[theorem]{Proposition}
\newtheorem{definition}[theorem]{Definition}

\newtheorem{conjecture}[theorem]{Conjecture}
\newtheorem{corollary}[theorem]{Corollary}

\newtheorem{remark}[theorem]{Remark}

\allowdisplaybreaks

\def\neweq#1{\begin{equation}\label{#1}}
\def\endeq{\end{equation}}
\numberwithin{equation}{section}
\begin{document}
{\centering
\bfseries
{\LARGE  Nash Williams Conjecture and the Dominating Cycle Conjecture\\ 
  }
 
\bigskip
\mdseries
{\large Arthur Hoffmann-Ostenhof}
\par
\upshape
\footnotesize \textit{Technical University of Vienna, Austria}\\

}

\noindent
\begin{abstract}
\noindent
The disproved Nash Williams conjecture states that every 4-regular 4-connected graph has a hamiltonian cycle. We show that a modification of this conjecture is equivalent to the Dominating Cycle Conjecture.
\end{abstract}

4-connected, cyclic 4-edge connected. 

\vspace{0.5cm}

\noindent

\section{Basic definitions and main result}

For used terminology which is not defined here we refer to \cite{Bo,BRV}.
 A \textit{dominating cycle} (DC) of a graph $G$ is a cycle which contains at least one endvertex of every edge of $G$. Let $v \in V(G)$ then $E_v$ denotes the set of edges incident with $v$. A \textit{closed trail} is a closed walk in which all the edges are distinct. All graphs here are considered to be loopless and finite.	

\noindent 
The following two conjectures are well known in graph theory. 
The first one was disproved by Meredith, see \cite {M}.


\noindent
\textit {Nash Williams Conjecture} (\textbf{NWC}): Every $4$-regular $4$-connected graph has a hamiltonian cycle.

\noindent
\textit {Dominating Cycle Conjecture} (\textbf{DCC}): Every cyclically $4$-edge connected cubic graph has a dominating cycle.

\noindent
The DCC is open and so far there is neither a promising approach known to prove it nor to disprove it.
For a survey on this conjecture, we refer to \cite{BRV}.

\noindent
We need the following definitions for introducing the modified NWC.

\begin{definition}\label{trans}
Let $H$ be a $4$-regular graph $H$ with a transition system $T$, i.e. $T:=\bigcup _{v \in V(H)} \{P_v\}$ where $P_v$ is a partition of 
the four edges incident with $v$ into two sets of size $2$; each of these two sets is called a \textbf{transition} of $T$, of $P_v$ and of $v$.
A trail is said to \textbf{follow a transition} if the two edges of the transition are consecutive edges of the trail.
Moreover, $H$ is said to be $\textbf{T-hamiltonian}$ if $H$ contains a \textbf{T-trail}, that is a spanning closed trail $C$ of $H$ such that for each $v \in V(H)$ one of the following two conditions is fulfilled:\\
a) $|E(C) \cap E_v|=2$ (in this case $C$ may follow no transition of $v$)\\
b) $C$ follows both transitions of $v$ (in this case $|E(C) \cap E_v|=4$). 
\end{definition}

\noindent
For an example see Figure 1. Observe that $H$ is $T$-hamiltonian if $H$ is hamiltonian. Hence, $T$-hamiltonicity generalizes the concept of hamiltonian graphs. Now, we introduce the modification of the NWC.

\begin{figure}[htpb]
\centering\epsfig{file=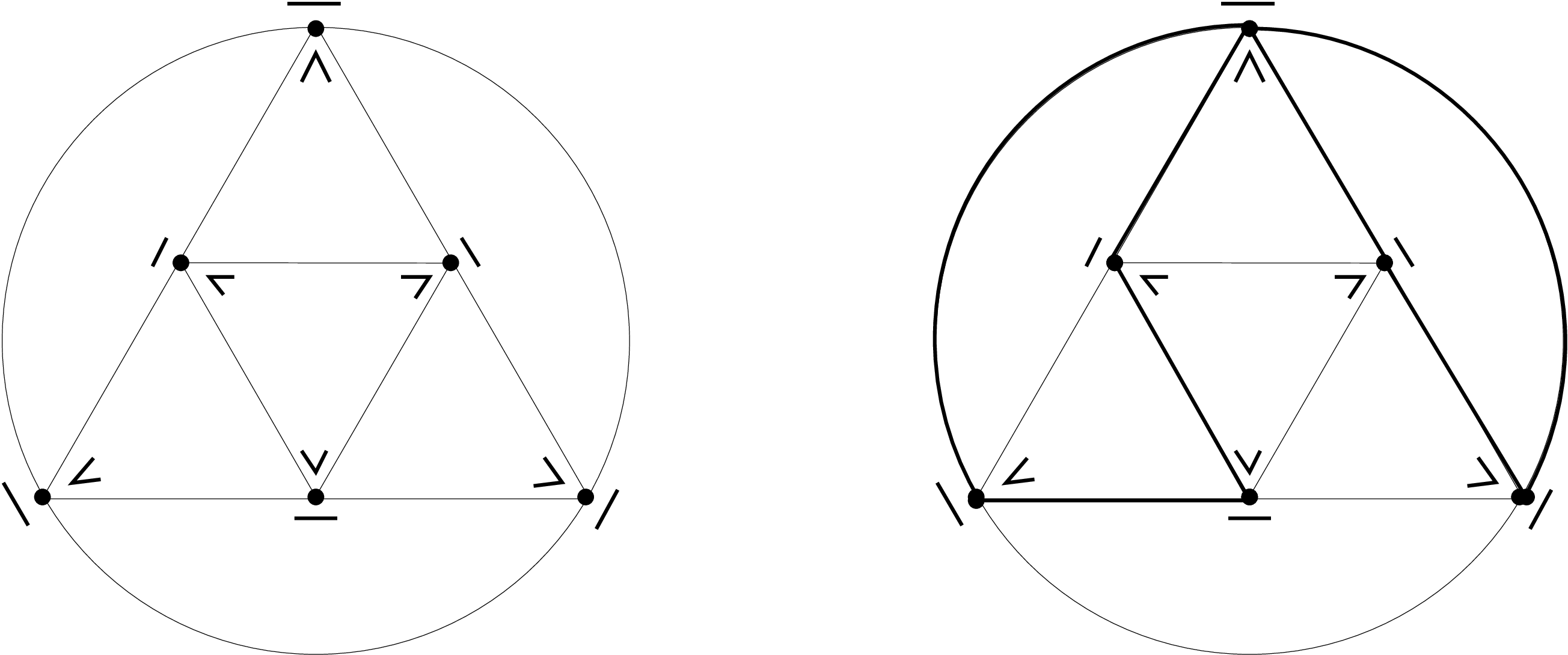,width=5.3 in}
\caption{On the left, a graph with a transition system (transitions are illustrated by short bold lines next to those pairs of edges which form a transition) and on the right, a T-trail (illustrated in bold) of the same graph with respect to its transition system.} 
\end{figure}

\noindent
\textit{\textbf{NWC*}}: Let $H$ be a $4$-regular $4$-connected graph $G$ with a transition system $T$, then $H$ is $T$-hamiltonian.  

\noindent
We state the main result.

\begin{theorem}\label{t1}
The DCC is equivalent to the NWC*.
\end{theorem}

\smallskip

\begin{remark}
Consider the NWC* as false. Then by the previous observation a counterexample to the NWC* is also a
counterexample to the NWC. Hence in order to disprove the NWC* one needs a counterexample to the NWC with additional properties. Therefore and by Theorem \ref{t1}, the DCC is harder to disprove than originally the NWC.  
\end{remark}

\noindent
Theorem \ref{t1} implies another result. For stating it we use the following definition.

\begin{definition}\label{dom}
We say that a cycle \textit{dominates} a matching $\mathcal M$ of a graph $G$ if the cycle contains at least one endvertex of 
each edge of $\mathcal M$. Let $G / \mathcal M$ denote the graph which results from $G$ by contracting each edge of $\mathcal M$ to a distinct vertex.
\end{definition} 

\noindent
By Corollary \ref{t2} of Theorem \ref{t1}, the subsequent conjecture is equivalent to the DCC.

\begin{conjecture}\label{con}
Let $G$ be a cubic graph with a perfect matching $\mathcal M$ such that $ G / \mathcal M$ is $4$-connected. Then $G$ contains a cycle which dominates $\mathcal M$.
\end{conjecture}

\section{Proof of the main result}


\noindent
A graph $G$ is called \textit{$k$-vertex connected} if $|V(G)| > k$ and $G-X$ is connected for every $X \subseteq V(G)$ with $|X|< k$. We abbreviate $k$-vertex connected by \textit{$k$-connected}. A graph $G$ is called \textit{$k$-edge connected} if $|V(G)| > 1$ and $G-Y$ is connected for every $Y \subseteq E(G)$ with $|Y|< k$. A set $E' \subseteq E(G)$ of a connected graph $G$ is called an \textit{edge cut} of $G$ if $G-E'$ is disconnected. Moreover, if at least two components of $G-E'$ are not a tree then $E'$ is also called a \textit{cyclic edge cut} of $G$.
If $G$ contains two vertex disjoint cycles, then $\lambda_c (G)$ is the minimum size over all cyclic edge cuts of $G$. The \textit{line graph} of a graph $G$ is denoted by $L(G)$. A cycle of length $3$ is called a \textit{triangle}.

\noindent
For convenience, we split the statement of Theorem \ref{t1} into Proposition \ref{p1} and Proposition \ref{p2}.

\begin{proposition}\label{p1}
If the NWC* is true, then the DCC is true.
\end{proposition}

\noindent
Proof. Let $G$ be a cubic graph satisfying $\lambda_c (G) \geq 4$. We show that $G$ has a DC. It is well known and not difficult to see that $L(G)$ is $4$-regular and $4$-connected. Note that every vertex $v \in V(G)$ corresponds to a unique triangle $t_v$ of $L(G)$. Moreover, $\{E(t_v)\,:\, v \in V(G)\}$ is a partition of $E(L(G))$. 

\noindent
Define the following transitions which imply a transition system $T$ of $L(G)$: each pair of edges which is incident with the same vertex of $L(G)$ and belongs to the same triangle $t_v$ of $L(G)$ forms a transition, see Figure 2.

\noindent
Since the NWC* holds, $L(G)$ is $T$-hamilitonian. If $L(G)$ has a hamiltonian cycle, then $G$ has a DC by Th. 5 in \cite{BRV} (see also \cite{HN}) and we are finished.

\begin{figure}[htpb]
\centering\epsfig{file=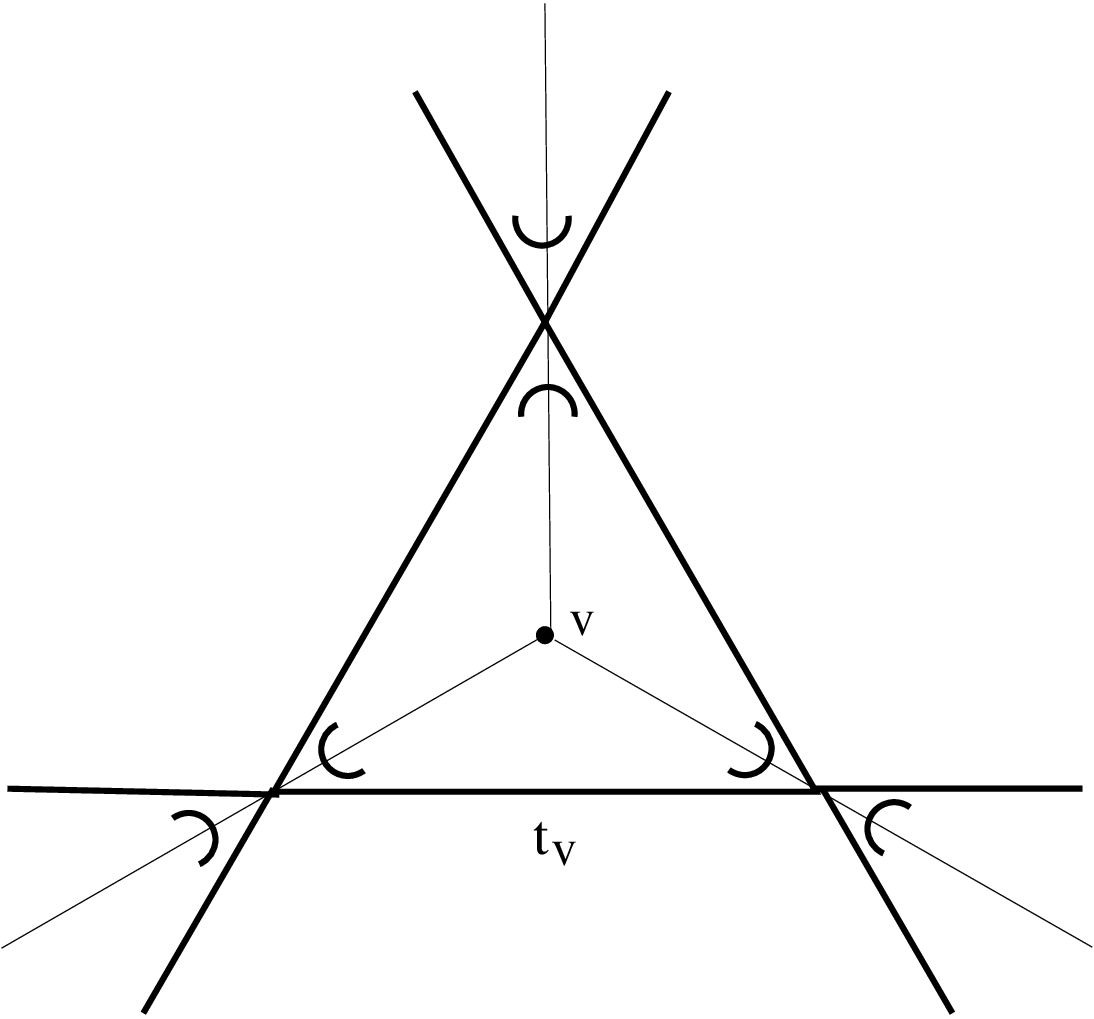,width=3.3 in}
\caption{The transitions of the three vertices of $t_v$ in the line graph.} 
\end{figure}

\noindent
Hence we may assume that $L(G)$ has a $T$-trail $C$ which contains a positive minimum number of $4$-valent vertices. $C$ contains at most one $4$-valent vertex in every triangle $t_v$ of $L(G)$. Otherwise $C$ follows by Def.\ref{trans} all transitions in at least two vertices of $t_v$ which implies that also the remaining third vertex of $t_v$ is a $4$-valent vertex of $C$. However this is not possible since $C$ obviously cannot follow now both transitions of the third vertex contradicting Def.\ref{trans}. Let $w$ be a $4$-valent vertex of $C$ and let $t_x$ with $x \in V(G)$ be one triangle of $L(G)$ which contains $w$. Define the new $T$-trail $C'$ of $L(G)$ which results from $C$ by replacing the two edges of $C$ which are contained in $t_x$ and incident with $w$ by the remaining edge of $t_x$. Then $C'$ contains fewer $4$-valent vertices than $C$ which contradicts the definition of $C$ and thus finishes the proof.


%

\noindent
For the formulation and proof of the next results we use the following definition.

\begin{definition}\label{def}
Let $H$ be a $4$-regular graph with a transition system $T$. Let $\mathcal G(H,T)$ denote the cubic graph which results from $H$ by firstly splitting every vertex $v$ of $H$ into two vertices $v'$, $v''$ of degree $2$ such that the corresponding edges of each transition of $T$ remain adjacent 
and by secondly adding the new edge $e(v):=v'v''$. Set $M_T:=\{ e(v) \,:\, v \in V(H)\}$ which is a perfect matching matching of $\mathcal G(H,T)$ satisfying $\mathcal G(H,T)/M_T \cong H$.
\end{definition}

\noindent
The following lemma can be verified straightforwardly.

\begin{lemma}\label{l3}
 Let $H$ be a $4$-regular graph with a transition system $T$. Then $\mathcal G(H,T)$ has a cycle which dominates $M_T$ if and only if 
$H$ has a $T$-trail.  
\end{lemma}

\noindent
We apply the subsequent three lemmata for the proof of the next proposition.



\begin{lemma}\label{known}
 Let $E_0$ be a matching of a cubic graph $G$ such that $G-E_0$ consists of two components $X_1$ and $X_2$.
 Moreover, let every edge of $E_0$ have precisely one endvertex in $X_i$, $i =1,2$. Then $|E_0 \cap M| + |E_0|$ is even for every perfect matching $M$ of $G$.
\end{lemma}

\noindent
Proof: straightforward. 

\begin{lemma}\label{l2}
Let $H$ be a $4$-regular $4$-connected graph with a transition system $T$, then $\mathcal G(H,T)$ is $3$-edge connected. Moreover, if $\mathcal G(H,T)$ has a cyclic $3$-edge cut $E_0$, then one of the two components of $\mathcal G(H,T)-E_0$ is a triangle. 
\end{lemma}

\noindent
Proof. Set $G:=\mathcal G(H,T)$. Obviously, $G$ is connected.
Suppose $G$ has a bridge $f$, then by Lemma \ref{known} (with $E_0 := \{f\}$), $f \in M_T$. Hence, $f$ corresponds to a cut vertex of $H$ which contradicts that $H$ is $4$-connected. 

\noindent
Suppose $E'$ is a $2$-edge cut of $G$. Then $E'$ is a matching, otherwise $G$ has a bridge. By Lemma \ref{known} we have two cases. 

\noindent
\textit{Case 1.} $E' \subseteq M_T$. \\
$H$ is a simple graph otherwise $H$ is not $4$-connected.
Thus, $G$ is also a simple graph. Therefore both components of $G-E'$ contain more than two vertices. Hence, $E'$ corresponds to a vertex $2$-cut of $H$ which contradicts the definition of $H$.
   
\noindent
\textit{Case 2.} $M_T \cap E' = \emptyset$.\\
Then $E'$ corresponds to a $2$-edge cut of $H$. Since $H$ is $4$-connected and since the edge connectivity is greater or equal the vertex connectivity, this is impossible.

\noindent
Hence, $G$ is $3$-edge connected which finishes the first part of the proof.

\noindent
Let $E_0$ be a cyclic $3$-edge cut of $G$. Then $E_0$ is a matching and $G-E_0$ consists of two components, otherwise $G$ would contain a $k$-edge cut for some $k \in \{1,2\}$. Moreover, at most one component of $G-E_0$ is a triangle since otherwise $|V(H)|=3$ and $H$ is not $4$-connected.

\noindent
Suppose by contradiction that no component of $G-E_0$ is a triangle. Hence both components have more than three vertices. Since every graph contains an even number of vertices of odd degree, every component has at least five vertices.  
Denote one of the two components of $G-E_0$ by $L$.
By Lemma \ref{known}, we need to consider two cases: Case A and Case B. Set $E_0 = \{e_1,e_2,e_3 \}$.

\noindent
\textit{Case A.} $M_T \cap E_0 = \{e_1\}$.

\noindent
Let us suppose first that $|L|>5$. \\
Set $E^*:= \{ \,e \in E(L) \cap M_T \,:\, e \,\, \text{is incident with a 2-valent vertex of L} \,\}$. Obviously, $1 \leq |E^*| \leq 2$ (depending on whether one edge of $M_T \cap E(L)$ covers one endvertex of $e_2$ and one of $e_3$).
Hence $E^* \cup e_1$ corresponds to a $j$-vertex cut of $H$ for some $j \in \{2,3\}$ which contradicts that $H$ is $4$-connected. 

\noindent
Thus, we may assume that $|V(L)| = 5$. Hence, $|M_T \cap E(L)| =2$. Set $\{a,b\}:= M_T \cap E(L)$. Since $H$ is a simple graph, $a$ and $b$ have the following properties:\\
\textbf{(1)} $a$ and $b$ are not contained together in a cycle of length $4$.\\
\textbf{(2)} Neither $a$ nor $b$ is contained in a triangle.

\noindent
Denote by $v_1$ the unique $2$-valent vertex of $L$ which is neither matched by $a$ nor by $b$. 
Denote by $a_1$ and $b_1$ the remaining two $2$-valent vertices in $L$. If $a_1b_1 \in E(L)$, then $a_1b_1 \not\in M_T$; otherwise the two vertices in $H$ corresponding to $e_1$ and $a_1b_1$ would form a vertex $2$-cut of $H$. 
Hence we can set $a:=a_1a_2$ and $b:= b_1b_2$ such that $V(L)= \{v_1,a_1,a_2,b_1,b_2\}$.

\noindent
By \textbf{(2)}, $v_1$ is adjacent to precisely one endvertex of $a$ and to one of $b$.
Thus, we have three cases:

\noindent
\textit {Case 1.} $\{v_1a_1, v_1b_1\} \subseteq E(L)$.\\
Then $G$ must contain the edge $a_2b_2$ twice which is impossible since $G$ is a simple cubic graph.

\noindent
\textit {Case 2.} $\{v_1a_1, v_1b_2\} \subseteq E(L)$.\\
Then $L$ contains the triangle consisting of the vertices $b_1$, $b_2$, $a_2$ which contradicts \textbf {(2)}.

\noindent
\textit {Case 3.} $\{v_1a_2, v_1b_2\} \subseteq E(L)$.\\
Then $L$ contains either double edges or a cycle of length $4$ consisting of the vertices: $a_1,a_2,b_1,b_2$ which contradicts \textbf{(1)}.

\noindent
Hence, Case A cannot occur.
 
\noindent
\textit{Case B.} $|M_T \cap E_0|=3$.
Since both components of $G-E_0$ have at least five vertices and since $E_0 \subseteq M_T$, $E_0$ corresponds to to a vertex $3$-cut of $H$ which 
contradicts the definition of $H$ and thus finishes the proof. 
 
\begin{definition}\label{dd}
Let $H$ be a $4$-regular simple graph with a transition system $T$.
Denote by $H'$ the graph which results from $\mathcal G(H,T)$ by contracting every triangle of $\mathcal G(H,T)$ to a distinct vertex. 
\end{definition}


\noindent
Note that each pair of triangles of $\mathcal G(H,T)$ (in Def.\ref{dd}) is vertex disjoint.
Hence, $H'$ is well defined.

\begin{lemma}\label{lll}
Let $H$ be a $4$-regular $4$-connected graph with a transition system $T$, then either 
$H' \in \{K_4, K_{3,3}\}$ or $\lambda_c (H') \geq 4$.
\end{lemma}

\noindent
Proof. Since $|V(H)| \geq 5$, $|V(\mathcal G(H,T))| \geq 10$ and thus $|V(H')| \geq 4$.
Suppose first that $H'$ does not contain two disjoint cycles. 
Since $|V(H')| \geq 4$ and by Theorem 1.2 in \cite{MC}, it follows that $H' \in \{K_4, K_{3,3}\}$.

\noindent
Now, assume that $H'$ has two disjoint cycles. Suppose by contradiction that $E'_0$ is a cyclic $k$-edge cut of $H'$ for some $k \in \{1,2,3\}$. Then $E'_0$ corresponds to a cyclic $k$-edge cut $E_0$ of $\mathcal G(H,T)$. Lemma \ref{l2} implies that $k=3$ and that one of the two components of $\mathcal G(H,T)-E_0$ is a triangle. Since this triangle is contracted to a vertex in $H'$, one of the two components of $H'-E'_0$ is a vertex. Hence, $E'_0$ is not a cyclic edge cut which is a contradiction and thus finishes the proof.

\begin{proposition}\label{p2}
If the DCC is true, then the NWC* is true.
\end{proposition}

\noindent
Proof. Let $H$ be a $4$-regular $4$-connected graph with a transition system $T$. Since $H$ is $4$-connected, $|V(H)| \geq 5$. Set $G:= \mathcal G(H,T)$. Then $|V(G)| \geq 10$ and we have the following two cases.

\noindent
\textit{Case 1.} $\lambda_c (G) \geq 4$. Then, by assumption $G$ has a DC which thus dominates $M_T$ (Def. \ref{def}). By Lemma \ref{l3}, $H$ is $T$-hamiltonian.

\noindent
\textit{Case 2.} $\lambda_c (G) < 4$. 
Consider $H'$. Every edge of $H'$ corresponds to an edge of $E(G)- \{e \in E(G)\, : \,e \,\,\text{is contained in a triangle of } G\}$.
Thus, every subgraph $X'$, say, of $H'$ induces by its corresponding edge set in $G$, a subgraph of $G$ which we denote by $X$. 

\noindent
Note that $K_4$ and $K_{3,3}$ have a dominating cycle. Therefore, and by Lemma \ref{lll} and since the DCC holds, $H'$ has a dominating cycle $C'$.
The corresponding subgraph $C \subseteq G$ is not a cycle if and only if there is a vertex $v' \in V(C')$ which has been obtained by contracting a triangle in
$G$. We denote this triangle by $\triangle(v')$ and call such a vertex $v'$, a \textit {bad vertex} of $C'$. We define the cycle $\widetilde C \subseteq G$ depending on $C'$: 
for each bad vertex $v'$ of $C'$, we extend $C \subseteq G$ to $\widetilde C$ by adding the unique path of length $2$ which is contained in $\triangle (v')$ and which connects two endvertices of two edges of $C$; if $C'$ has no bad vertex, then $C$ is already a cycle and we set $\widetilde C:= C$. We show that $\widetilde C$ dominates $M_T$.  

\noindent
Since $H$ is $4$-connected, $H$ is simple. Thus, no triangle of $G$ contains an edge of $M_T$. Hence, it suffices to show that the edge $e \in E(G)$ is dominated by $\widetilde C $ for every $e' \in E(H')$. If $e' \in E(C')$, then $e \in E(C)$ and thus $e$ is dominated by $\widetilde  C$. If $e' \not\in E(C')$, then two cases are possible.

\noindent
\textit{Case A.} $e'$ is a chord of $C'$.\\ Then, by the construction of $\widetilde C$ both endvertices of $e$ are contained in $\widetilde C$.

\noindent
\textit{Case B.} Precisely one endvertex of $e'$ is contained in $C'$.\\ Then, one endvertex of $e$ is contained in $\widetilde C$.

\noindent
Hence, $\widetilde C$ dominates $M_T$. By Lemma \ref{l3}, $H$ is $T$-hamiltonian which finishes the proof.


\noindent
\textbf{Theorem \ref {t1}}  \textit {The DCC is equivalent to the NWC*.}
 
\smallskip

\noindent
Proof. By Prop. \ref{p1} and Prop. \ref{p2}, Theorem \ref{t1} follows.




\begin{corollary}\label{t2}
Conjecture \ref{con} is equivalent to the DCC.
\end{corollary}

\noindent
Proof. Suppose that Conjecture \ref{con} holds. We show that in this case the NWC* holds which implies by Th.\ref{t1} the truth of the DCC.
Let $H$ be a $4$-regular $4$-connected graph with a transition system $T$. 
Since Conjecture \ref{con} holds, the cubic graph $\mathcal G(H,T)$ has a cycle dominating $M_T$ and thus by Lemma \ref{l3}, the first part of the proof is finished. 

\noindent
Suppose that the DCC and thus by Th.\ref{t1} also the NWC* holds. Let $G$ and $\mathcal M$ be defined as in Conjecture \ref{con}. We want to find a cycle which dominates $\mathcal M$. Set $H:= G/\mathcal M$ and let $T$ be the transition system of $H$ such that $\mathcal G(H,T) \cong G$.
Since the NWC* holds, $H$ is $T$-hamiltonian. By Lemma \ref{l3}, $G$ has a cycle dominating $M_T$ where $M_T$ equals $\mathcal M$ which finishes the proof.

\smallskip
 
\noindent
{\textbf{Acknowledgment}: This work was funded by the Austrian Science Fund (FWF-Project 26686).

\footnotesize
\bibliographystyle{plain}

\end{document}